\documentclass[12pt]{article}
\usepackage{latexsym, amsmath, amssymb, amsthm}
\renewcommand\le{\leqslant}
\renewcommand\ge{\geqslant}
\newtheorem{theorem}{Theorem}
\newtheorem{lemma}{Lemma}
\renewenvironment{proof}{\par\medskip\noindent{\bf Proof.}}%
{\nopagebreak\hfill$\Box$\par\medskip}

\begin{document}
\title{Nonantagonistic noisy duels of discrete type \\
with an arbitrary number of actions}
\author{Lyubov~N.~Positselskaya}
\date{{\it  Moscow State Socio-Humanitary Institute\\
e-mail: posicelskaja@do.tochka.ru}}
\maketitle

\begin{abstract}
We study a nonzero-sum game of two players which is
a generalization of the antagonistic noisy duel of discrete type.
The game is considered from the point of view of various criterions
of optimality.
We prove existence of $\varepsilon$\nobreakdash-equilibrium situations
and show that the $\varepsilon$\nobreakdash-equilibrium strategies
that we have found are $\varepsilon$\nobreakdash-maxmin.
Conditions under which the equilibrium plays are Pareto-optimal
are given.
\end{abstract}

\section {Introduction}
   The classical duel is a zero-sum game of two players
of the following type.  The players have certain resources
and use them during a given time interval with the goal of
achieving success.  Use of the resource $\gamma$ at the moment~$t$
leads to success with the probability depending on the amount of
resource $\gamma$ and the time $t$ only (it is usually assumed that
the probability of success increases with time).  As soon as one
player achieves the goal he receives his profit, which is equal to
his opponent's loss, and the game ends.
   Various assumptions about the ways players use their resources and
about players receiving information about the opponent's behavior
during the game define various kinds of
duels~\cite{kar,kimoverview}.
   Models were considered where players' resources were discrete
({\it discrete firing\/} duels), infinitely divisible
({\it continuous firing\/} duels), continuous for one of the players
and discrete for the other one ({\it mixed\/} duels, or
{\it fighter-bomber duels\/})~\cite{pdn,pnr}.
   Researchers studied {\it noisy\/} duels~\cite{fk,pdn,ran},
where every player at a given moment of time had complete information
about his opponent's behavior up to that moment, and {\it silent\/}
duels, where no such information was available.
   At present time duels are considered as classical models of
competition~\cite{kar,pmixsin}.  However, their application as such
is somewhat limited by the assumption that the players' interests
are strictly opposite to each other.
   Games of the duel type with nonzero sum belong to an unexplored
class of infinite games with nonopposite interests.

   A nonzero-sum game which is a generalization of the classical
antagonistic duel was first considered in~\cite{pnozmixn}; it was
a nonantagonistic noisy fighter-bomber duel.  Then nonzero-sum duels
were studed in~\cite{ndst},~\cite{fpm},~\cite{pareto}.

   In this paper we study a nonzero-sum game of two players
which is a generalization of the antagonistic noisy duel with
discrete firing~\cite{fk,fox}.
   This paper is an extended verion of~\cite{fpm2} with
complete proofs.
   It is also a sequel to~\cite{fpm}.

   The author is grateful to Leonid Positselski for his help in
translating this paper into English and editing it.

\section {Preliminaries from Game Theory}
 A {\it game of two players\/} is a quadruple
$$\Gamma=\{X,Y,K_1(x,y), K_2(x,y)\},$$
where $X$ and $Y$ are sets of the players' {\it strategies\/}
and $K_j(x,y),\;j=1,2$ are the players' {\it payoff functions},
which are defined on the Carthesian product $X\times Y$ and
determine the $j$-th player's payoff when Player~I uses
a strategy $x\in X$ and Player~II uses a strategy $y\in Y$.
 A game $\Gamma$ is called a {\it zero-sum game of two players},
or an {\it antagonistic game\/} if $K_1(x,y)+K_2(x,y)=0$
for all $x\in X$, $y\in Y$ and a {\it nonzero-sum game of
two players}, or a {\it nonantagonistic game\/} otherwise.
 The {\it mixed extension\/} of a game $\Gamma$ is the game
$\overline\Gamma=\left(\Phi,\Psi, \overline{K}(\varphi,\psi)\right)$,
where $\Phi$ and $\Psi$ are the sets of distributions on $X$ and $Y$
and $\overline{K}_j(\varphi,\psi)$ ($j=1,2$) are the mean values
of the payoff functions $K_j(x,y)$ over the distributions
$\varphi\in\Phi$, $\psi\in\Psi$.
 {\it Pure strategies\/} of the game $\overline\Gamma$ are distributions
concentrated in one point $x\in X$ or $y\in Y$, respectively.
All the other distributions $\varphi\in\Phi$, $\psi\in\Psi$ are called
{\it mixed strategies}.

A {\it situation\/} $(x,y)$ is a pair of the player's strategies.
  A situation  $(x_e,y_e)$ is called an {\it equilibrium situation\/}
if for all $x\in X$, $y\in Y$ the following inequalities hold:
\begin{align}
&K_1(x,y_e)\leqslant K_1(x_e,y_e);\qquad K_2(x_e,y)
\leqslant K_2(x_e,y_e).\label {1}
\end{align}
The vector $(v^1, v^2)$, where $v^1=K_1(x_e,y_e)$ and $v^2=K_2(x_e,y_e)$,
is called the {\it equilibrium value\/} corresponding to the equilibrium
situation $(x_e,y_e)$, and a strategy included in some equilibrium
situation is called an {\it equilibrium\/} strategy.

An infinite game with nonantagonistic interests may not admit
an equilibrium situation.
  A situation $(x^{\varepsilon},y^{\varepsilon})$ is called
an {\it $\varepsilon$-equilibrium situation\/} if for any
$x\in X$, $y\in Y$ the following inequalities hold:
\begin{equation}
K_1(x,y^{\varepsilon})-\varepsilon\leqslant
K_1(x^{\varepsilon},y^{\varepsilon}); \qquad
K_2(x^{\varepsilon},y)-\varepsilon\leqslant
K_2(x^{\varepsilon},y^{\varepsilon}).
\label {3}
\end{equation}
A player's strategy included in some $\varepsilon$-equilibrium situation
will be called an {\it $\varepsilon$-equilibrium\/} strategy.
If the limits
$$
 v^j=\lim_{\varepsilon\to 0}K_j(x^{\varepsilon},y^{\varepsilon})
 \quad(j=1,2),
$$
exist, the vector $(v^1, v^2)$ will be called the {\it equilibrium
value\/} corresponding to a set of $\varepsilon$-equilibrium situations
$\{(x^{\varepsilon},y^{\varepsilon})\}$.

A strategy $x_m\in X$ is called a {\it maxmin\/} strategy of Player I
if the function $\displaystyle m_1(x)=\inf_{y\in Y}K_1(x,y)$ attains
its maximal value in it.
Analogously, a strategy $y_m\in Y$ is called a {\it maxmin\/} strategy
of Player II if the function
$\displaystyle m_2(y)=\inf_{x\in X}K_2(x,y)$
attains its maximal value in it.
The vector $$w=(w^1, w^2), \mbox { where }
 w^1=\max_{x\in X}\inf_{y\in Y}K_1(x,y), \
 w^2=\max_{y\in Y}\inf_{x\in X}K_2(x,y),
$$
is called the {\it maxmin value\/} of a nonzero-sum game.
The value $w^j$ is the best guaranteed payoff of the $j$-th player.

In an infinite game maxmin strategies may not exist. In this case
{\it $\varepsilon$-maxmin\/} strategies are considered.
A strategy $x_m^{\varepsilon}\in X$ is called
an {\it $\varepsilon$-maxmin\/} strategy of Player I if
$$m_1(x_m^{\varepsilon})>w^1-\varepsilon,  \text{ where }
  w^1=\sup_{x\in X}\inf_{y\in Y}K_1(x,y),
 $$
and $w^1$ is called the maxmin value of the game of Player I.
Analogously one defines an {\it $\varepsilon$-maxmin\/} strategy
of Player II:
a strategy $y_m^{\varepsilon}\in X$ is called
an {\it $\varepsilon$-maxmin\/} strategy of Player II if
$$m_2(y_m^{\varepsilon})>w^2-\varepsilon,  \text{ where }
 w^2=\sup_{y\in Y}\inf_{x\in X}K_2(x,y),$$
and $w^2$ is called the maxmin value of the game of Player II.

If an antagonistic game admits equilibrium situations, then
the equilibrium values corresponding to them coincide and are equal to
the maxmin value.
Equilibrium strategies of an antagonistic game are maxmin and they
are called {\it optimal}, while $\varepsilon$-equilibrium values are
$\varepsilon$-maxmin and they are called {\it $\varepsilon$-optimal}.
A nonzero-sum game may admit equilibrium situations with unequal
equilibrium values.

Let us introduce a partial order relation $\succ$ in $\mathbb R^2$:
$$a\succ b, \text{ if } a_j\geqslant b_j\; (j=1,2),$$
where at least one of the inequalities is strict.
Denote by $S$ the set of all situations $s=(x,y)$.
 The relation of {\it Pareto preference\/} on the set $S$ is defined
as follows.  Let $s^i=(x^i,y^i)$ ($i=1,2$).  Then
\begin{align}
 &s^1\succ s^2, \text{ if }
K^1\succ K^2, \text{ where }
K^i=(K_1(x^i,y^i), K_2(x^i,y^i)), \ i=1,2.
\end{align}
 A situation $s^p=(x_p,y_p)$ is called {\it Pareto-optimal\/} if
there are no situations $s=(x,y)$ such that $s\succ s^p$.

In an antagonistic game all the situations are Pareto-optimal.
A game is called {\it quasi-antagonistic}, or a {\it game with
opposite interests}, if all the situations of the game are
Pareto-optimal, i.~e., if for any two situations
$(x^1, y^1), (x^2, y^2)\in S$ the following conditions hold:
\begin{align}
&K_1(x^1, y^1)< K_1(x^2, y^2)\iff K_2(x^1, y^1)> K_1(x^2, y^2));
\nonumber \\
&K_1(x^1, y^1)= K_1(x^2, y^2)\iff K_2(x^1, y^1)= K_1(x^2, y^2)).
                                   \label {2a}
\end{align}

\section {Posing the problem}
 Consider a nonzero-sum game of two players which is
a generalization of the antagonistic noisy duel with discrete firing.

The players compete in the conditions of complete information.
They have discrete resources $m, \;n \in \mathbb N$ which they use
during the time interval $[0,1]$.
The effectiveness of the $j$-th player using his resource is described
by the function $P_j(t)$  ($j=1,2$), which defines the probability of
achieving success when using the unit of resource at the moment~$t$.
The functions $P_j(t)$ are continuous and increasing,
$P_j(0)=0$, $P_j(1)=1$, $0<P_j(t)<1$ for $t\in (0,1)$.
If one of the players achieves success, the game stops.
If a player has used all of his resource and hasn't achieved success,
the other player postpones his action until the moment $t=1$, when
his probability of success is equal to one.
The profit of the $j$-th player in the case of his success is equal to
$A_j$, and his loss in the case of his opponent's success is equal to
$B_j$, where
\begin{align}
 A_j\geqslant 0,\quad B_j\geqslant 0,\quad A_j+B_j>0,\quad j=1,2.
\label{I1}
\end{align}
The players' profits are equal to $0$ if no one of them achieved success
or if the success was achieved by both of them simultaneously.
A player's strategy is a function assigning the moment of next action to
a pair of amounts of players' current resources.
Let us call the described game a {\it noisy nonzero-sum duel
with discrete firing}.

Denote by $\tau_i$
($0\leqslant\tau_i\leqslant\tau_{i-1}\leqslant 1$, $i=1,2,\dots,m$)
the moments of time when Player I uses his resource.
Analogously, denote by $\eta_i$
($0\leqslant\eta_i\leqslant\eta_{i-1}\leqslant 1$, $i=1,2,\dots,n$)
the moments of time when Player II uses his resource.
The vectors $\tau$ and $\eta$ will be called the {\it vectors of
action moments\/} of the players.
The payoff function $K_j(\tau,\eta)$ is the mathematical expectation
of profit received by the $j$-th player in the case when the players
use their resources at the moments of time
$\tau_i$, $\eta_j$ ($i=1,\dots,m$; $j=1,\dots,n$).
It is computed in the following way.  If $m=0$, $n=0$, then $K_1=K_2=0$.
If $m\ge 1$, $n=0$, then $K_1=A_1$, $K_2=-B_2$.
If $m=0$, $n\ge 1$, then $K_1=-B_1$, $K_2=A_2$.
Assume that $m\ge 1$, $n\ge 1$. Set
$$\tau'=(\tau_{m-1},\dots,\tau_1);\quad
\eta'=(\eta_{n-1},\dots,\eta_1).$$
Then
\begin{equation}
K_1(\tau, \eta)=
 \begin{cases}
A_1P_1(\tau_m)+(1-P_1(\tau_m))K_1(\tau',\eta),
&\text{if }\tau_m<\eta_n,\\
A_1P_1(\tau_m)(1-P_2(\tau_m))-B_1(1-P_1(\tau_m))P_2(\tau_m)+\\
 +(1-P_1(\tau_m))(1-P_2(\tau_m))K_1(\tau',\eta'),
&\text{if }\tau_m=\eta_n,\\
-B_1P_2(\eta_n)+(1-P_2(\eta_n))K_1(\tau,\eta'),
&\text{if }\tau_m>\eta_n.
  \end{cases}
\label {4a}
\end{equation}

\begin{equation}
K_2(\tau; \eta)=
 \begin{cases}
A_2P_2(\eta_n)+(1-P_2(\eta_n))K_2(\tau, \eta'),
&\text{if }\eta_n<\tau_m,\\
A_2P_2(\tau_m)(1-P_2(\tau_m))-B_2(1-P_1(\tau_m))P_2(\tau_m)+\\
 +(1-P_1(\tau_m))(1-P_2(\tau_m))K_2(\tau',\eta'),
&\text{if }\tau_m=\eta_n,\\
-B_2P_1(\tau_m)+(1-P_2(\tau_m))K_2(\tau',\eta),
&\text{if }\eta_n>\tau_m.
  \end{cases}
\label {4}
\end{equation}

Let us denote $A=(A_1, A_2)$, $B=(B_1, B_2)$
and call $A$ the profit vector and $B$ the loss vector of the players.
Introduce the vector-function of effectiveness $P(t)=(P_1(t), P_2(t))$.
Let us denote the described game by $\Gamma_{mn}(P,A,B)$ or,
for brevity, by $\Gamma_{mn}$ and its mixed extension by
$\overline{\Gamma}_{mn}(P,A,B)$ (or $\overline{\Gamma}_{mn}$).

Suppose that in the duel $\Gamma_{mn}(P,A,B)$ the profit and loss
vectors are related by the equation $A_1=B_2$,  $A_2=B_1$, i.~e.,
the profit of each player is equal to his opponents' loss.
Then it follows from the relations (\ref{4a}), (\ref{4}) that
$K_1(\tau, \eta)=-K_2(\tau, \eta)$, i.~e., under these conditions
the game is antagonistic.

\section {Situations of $\protect\varepsilon$-equilibrium}

\begin{lemma}[Fox, Kimeldorf~\cite{fk}] \label{l.0}
 There exists a set $\{t_{ij}\; i,j\in \mathbb N\}$ such that
\begin{equation}
\prod_{i=1}^m\left(1-P_1(t_{in})\right)+
\prod_{j=1}^n\left(1-P_2(t_{mj})\right)=1,
\label {5}
\end{equation}
and for all $m,n\in\mathbb N$ the following inequalities hold:
 $$0<t_{mn}<\min(t_{m-1,n},t_{m,n-1}), \text{ where } t_{0n}=t_{m0}=1.$$
\end{lemma}
 Set
\begin{equation}
\lambda=\min\left\{1/(A_1+B_1),1/(A_2+B_2)\right\}/2.
\label {d.1}
\end{equation}
 Choose $\varepsilon>0$ and find $\delta_j$ ($j=1,2$) such that
\begin{align}
t_{mn}<\delta_j<\min(t_{m-1,n},t_{m,n-1}) \;\text{ and }
P_j(\delta_j)<P_j(t_{mn})+\lambda\varepsilon.
\label {d.2}
\end{align}
Take $\delta=\min\{\delta_1,\delta_2\}$.
Let $\mu$, $\nu$ be the current values of the players' resources
and $\varphi_{\mu\nu}^{\varepsilon}$ be the uniform distribution
concentrated in the interval $[t_{\mu\nu},t_{\mu\nu}+\delta]$.
Define the players' mixed strategies $x^{\varepsilon}$, 
$y^{\varepsilon}$ in the following way.
The strategy $x^{\varepsilon}$ ($y^{\varepsilon}$) prescribes
to choose the next moment of action $\tau_{\mu}$ ($\eta_{\nu}$)
in the random way according to the distribution function
$\varphi_{\mu\nu}^{\varepsilon}$.

\begin{theorem} \label{t.1}
The situations $(x^{\varepsilon}, y^{\varepsilon})$ in the game
$\overline{\Gamma}_{mn}(P,A,B)$ are $\varepsilon$-equilibrium situations.
The vector $v_{mn}=(v_{mn}^1, v_{mn}^2)$ defined by the formulas
\begin{align}
&v_{mn}^1=
A_1-(A_1+B_1)\prod_{i=1}^m\left(1-P_1(t_{in})\right)=
(A_1+B_1)\prod_{j=1}^n\left(1-P_2(t_{mj})\right)-B_1;
\nonumber
\\
&v_{mn}^2=
(A_2+B_2)\prod_{i=1}^m\left(1-P_1(t_{in}\right)-B_2=
A_2-(A_2+B_2)\prod_{j=1}^n\left(1-P_2(t_{mj})\right)
\label {6}
  \end{align}
is the corresponding equilibrium value.
\end{theorem}

\begin{proof}
 The second equations in both lines of~\eqref{6} hold by Lemma~1.
 Let us obtain recurrence relations for $v_{mn}$. By~\eqref{6} we have:
\begin{align*}
&v_{m-1,n}^1=
A_1-(A_1+B_1)\prod_{i=1}^{m-1}\left(1-P_1(t_{in})\right);
\\
&v_{mn}^1=
A_1-(A_1+B_1)\left(1-P_1(t_{mn})\right)
\prod_{i=1}^{m-1}\left(1-P_1(t_{in})\right)=
\\
&\phantom{v_{mn}^1}=A_1-\left(1-P_1(t_{mn})\right)(v_{m-1,n}^1-A_1).
  \end{align*}
Finally we have an expression for $v_{mn}^1$ in terms of $v_{m-1,n}^1$:
\begin{align}
&v_{0n}^1=-B_1\ \text { for } n>0;
\nonumber
\\
&v_{mn}^1=
A_1P_1(t_{mn})+\left(1-P_1(t_{mn})\right)v_{m-1,n}^1
\ \text { for } m>0,\;n>0.
\label {6a}
  \end{align}
Now obtain an expression for $v_{mn}^1$ in terms of $v_{m,n-1}^1$.
By~\eqref{6} we have:
\begin{align*}
&v_{m,n-1}^1=
(A_1+B_1)\prod_{j=1}^{n-1}\left(1-P_2(t_{mj})\right)-B_1;
\\
&v_{mn}^1
=(A_1+B_1)\left(1-P_2(t_{mn})\right)
\prod_{j=1}^{n-1}\left(1-P_2(t_{mj})\right)-B_1=
\\
&\phantom{v_{mn}^1}=\left(1-P_2(t_{mn})\right)(v_{m,n-1}^1+B_1)-B_1.
  \end{align*}
Hence we get:
\begin{align}
&v_{m0}^1=A_1\ \text{ for }m>0;
\nonumber
\\
&v_{mn}^1=-B_1P_2(t_{mn}+\left(1-P_2(t_{mn})\right)v_{m,n-1}^1
\ \text{ for } m>0,\;n>0.
\label {6b}
  \end{align}
Analogously we deduce recurrence relations for
$v_{mn}^2$ in terms of $v_{m-1,n}^2$ and $v_{m,n-1}^2$:
\begin{align}
&v_{m0}^2=-B_2\ \text{ for }m>0;
\nonumber
\\
&v_{mn}^2=
A_2P_2(t_{mn})+\left(1-P_2(t_{mn})\right)v_{m,n-1}^2
\ \text { for } m>0,\;n>0.
\label {6c}
\\
&v_{0n}^2=A_2\ \text{ for }n>0;
\nonumber
\\
&v_{mn}^2=-B_2P_1(t_{mn}+\left(1-P_1(t_{mn})\right)v_{m-1,n}^2
\ \text { for } m>0,\;n>0.
\label {6d}
  \end{align}

To prove that
 $(x^{\varepsilon}, y^{\varepsilon})$
are $\varepsilon$-equilibrium situations it is necessary and sufficient
to verify the following relations:
\begin{align}
&\overline{K}_1(\tau, y^{\varepsilon})\leqslant
\overline{K}_1(x^{\varepsilon}, y^{\varepsilon})+\varepsilon
\mbox{ for any } \tau;
\label {t1.1}
\\
&\overline{K}_2(x^{\varepsilon},\eta)\leqslant
\overline{K}_2(x^{\varepsilon},y^{\varepsilon})+\varepsilon
\mbox{ for any } \eta;
\label {t1.2}
\\
&\lim_{n\to \infty} \overline{K}_1(x^{\varepsilon}, y^{\varepsilon})
=v_{mn}^1;
\label {t1.3}
\\
&\lim_{n\to \infty} \overline{K}_2(x^{\varepsilon}, y^{\varepsilon})
=v_{mn}^2.
\label {t1.4}
\end{align}

\begin{lemma}\label{l.t1}
For any $\varepsilon>0$ there exist such strategies
$x^{\varepsilon}$, $y^{\varepsilon}$ that
for any pure strategies $\tau$ and $\eta$ of Players I and II
the following inequalities hold:
\begin{align}
&\overline{K}_1(\tau,y^{\varepsilon})<v_{mn}^1+\varepsilon;
\label {t1.5}
\\
&\overline{K}_1(x^{\varepsilon},\eta)> v_{mn}^1-\varepsilon.
\label {t1.6}
\end{align}
\end{lemma}

\begin{proof}
For $n=0$ and arbitrary $m>0$ by the definition of payoff function
we have:
$$K_1=A_1;\quad K_2=-B_2.$$
Analogously  for $m=0$ and arbitrary  $n>0$
$$K_1=-B_1;\quad K_2=A_2.$$
By~\eqref{6} we have:
$$v_{m0}^1=A_1;\quad v_{m0}^2=-B_2;\quad
v_{0n}^1=-B_1;\quad v_{0n}^2=A_2.$$
In both cases the inequalities~\eqref{t1.5}, \eqref{t1.6} hold.

For arbitrary $m>0$, $n>0$ we proceed by induction on the number of
action moments of the players.
Assume that Lemma is true for all pairs $(\mu,\nu)$ for which
$\mu\le m$, $\nu\le n$, $(\mu,\nu)\neq(m,n)$ and prove it for
$(\mu,\nu)=(m,n)$.

Every pure strategy of Player~I has the following structure.
Let $t\in[0,1]$ denote the planned moment of his first action.
If Player~II acts and misses at a time $\eta_n<t$, then Player~I
follows a pure strategy $\tau^1$ in $\Gamma_{m,n-1}$.
If Player~II does not act before the time $t$ then Player~I acts at
the time $t$, and unless Player~II also acts at the time $t$,
Player I after that adopts a pure strategy $\tau^2$ in $\Gamma_{m-1,n}$.

Every pure strategy of Player~II has the similar structure.
Let $u\in[0,1]$ denote the planned moment of his first action.
If Player~I acts and misses at a time $\tau_m<u$, then Player~II
follows a pure strategy $\eta^1$ in $\Gamma_{m-1,n}$. If Player~I
does not act before the time $u$ then Player~II acts at the time $u$,
and unless Player~I also acts at the time $u$, Player~II after that
adopts a pure strategy $\eta^2$ in $\Gamma_{m,n-1}$.

The strategy $y^{\varepsilon}$ of Player II is constructed as follows.
Fix $\varepsilon>0$.  Choose $t$ randomly according to
$\varphi_{mn}^{\varepsilon}$.  If Player I does not act before
the time $t$, then Player II acts at the time $\eta_n=t$ and, unless
Player I also acts at the time $t$, then adopts the strategy
$y^{1,\varepsilon}$ in $\Gamma_{m,n-1}$.
 If Player I acts and misses at the time $\tau_m<u$, then Player II
adopts the strategy $y^{2,\varepsilon}$ in $\Gamma_{m-1,n}$.
Accoding to the inductive assumption, we choose $y^{1,\varepsilon}$,
$y^{2,\varepsilon}$ such that for any pure strategies $\tau^1$ and
$\tau^2$ of Player I the following inequalities hold:
\begin{align}
&\overline{K}_1(\tau^1,y^{1,\varepsilon})<v_{m,n-1}^1+\varepsilon/2
\label {lt.3}
\\
&\overline{K}_1(\tau^2,y^{2,\varepsilon})<v_{m-1,n}^1+\varepsilon/2
\label {lt.4}
\end{align}
 The strategy $x^{\varepsilon}$ of Player I is constructed similarly.
Choose $u$ randomly according to $\varphi_{mn}^{\varepsilon}$.
If Player II does not act before the time $u$, then Player II acts at
the time $\eta_n=u$ and, unless Player II also acts at the time $u$,
adopts the strategy $x^{1,\varepsilon}$ in $\Gamma_{m-1,n}$.
 If Player II acts and misses at the time $\tau_m<u$, then Player I
adopts the strategy $x^{2,\varepsilon}$ in $\Gamma_{m,n-1}$.
Accoding to the inductive assumption, we choose  $x^{1,\varepsilon}$,
$x^{2,\varepsilon}$ such that for any pure strategies $\eta^1$ and
$\eta^2$ of Player II the following inequalities hold:
\begin{align}
&\overline{K}_1(x^{1,\varepsilon},\eta^1)> v_{m-1,n}^1-\varepsilon/2
\label {lt.5}
\\
&\overline{K}_1(x^{2,\varepsilon},\eta^2)> v_{m,n-1}^1-\varepsilon/2
\label {lt.6}
\end{align}
For all strategies described above we have ignored the response to
simultaneous actions of players since this event has probability 0
when Player I adopts $x^{\varepsilon}$
or Player II adopts $y^{\varepsilon}$.

Let us prove the inequality~\eqref{t1.5}.

Note that by the definition of payoff function for any strategies
$x$, $y$ we have:
\begin{align}
   &-B_j\le\overline{K}_j(x,y)\le A_j\;(j=1,2)\implies
\nonumber
 \\
   &\overline{K}_j(x,y)+B_j\ge 0;
\label{lt.1}
\\
   &A_j-\overline{K}_j(x,y)\ge 0.
\label{lt.2}
\end{align}

Let $\tau$ be an arbitrary pure strategy of Player I and $t$ be
the time of first action in $\tau$.
There are three cases to be considered.
\begin{enumerate}
\item
Suppose $t\in [0,t_{mn}]$.
In this case Player I always acts first; so
\begin{align}
\overline{K}_1(\tau,y^{\varepsilon})&=
A_1P_1(t)+(1-P_1(t))\overline{K}_1(\tau^2,y^{2,\varepsilon}).
\label{lt.2a}
  \end{align}
One can see from the inequality~\eqref{lt.2} that the right hand side
of~\eqref{lt.2a} does not exceed
\begin{align*}
&(A_1-\overline{K}_1(\tau^2,y^{2,\varepsilon}))P_1(t_{mn})+
\overline{K}_1(\tau^2,y^{2,\varepsilon})=
\\
&=A_1P_1(t_{mn})+(1-P_1(t_{mn}))\overline{K}_1(\tau^2,y^{2,\varepsilon}).
  \end{align*}
According to~\eqref{lt.4}, the latter expression is less than
\begin{align*}
 A_1P_1(t_{mn})+(1-P_1(t_{mn}))(v_{m-1,n}+\varepsilon/2).
  \end{align*}
Using the recurrence relation~\eqref{6a} we get:
\begin{align*}
 &A_1P_1(t_{mn})+(1-P_1(t_{mn}))(v_{m-1,n}+\varepsilon/2)=
v_{mn}+(1-P_1(t_{mn}))\varepsilon.
  \end{align*}
Therefore
\begin{align*}
\overline{K}_1(\tau,y^{\varepsilon})<v_{mn}+\varepsilon.
  \end{align*}
\item
Suppose $t\in (t_{mn},t_{mn}+\delta)$.
In this case either player may act first; so
\begin{align*}
\overline{K}_1(\tau,y^{\varepsilon})&=
\int\limits_{t_{mn}}^t(-B_1P_2(\xi)+(1-P_2(\xi))\overline{K}_1
(\tau^1,y^{1,\varepsilon}) )d\varphi_{mn}^{\varepsilon}(\xi)+
\\
&+\int\limits_t^{t_{mn}+\delta}
(A_1P_1(\xi)+(1-P_1(\xi))\overline{K}_1(\tau^2,y^{2,\varepsilon}) )
d\varphi_{mn}^{\varepsilon}(\xi).
  \end{align*}
Estimate the first integrand using the relations~\eqref{lt.1},
 \eqref{lt.3}, and \eqref{6b}:
\begin{align*}
&-B_1P_2(\xi)+(1-P_2(\xi))\overline{K}_1(\tau^1,y^{1,\varepsilon})=
\\
&=\overline{K}_1(\tau^1,y^{1,\varepsilon})- 
P_2(\xi)(\overline{K}_1(\tau^1,y^{1,\varepsilon})+B_1)\le
\\
&\le\overline{K}_1(\tau^1,y^{1,\varepsilon})-
P_2(t_{mn})(\overline{K}_1(\tau^1,y^{1,\varepsilon})+B_1)=
\\
&-B_1P_2(t_{mn})+(1-P_2(t_{mn}))
\overline{K}_1(\tau^1,y^{1,\varepsilon})\le
 \\
&\le -B_1P_2(t_{mn})+(1-P_2(t_{mn}))(v_{m,n-1}+\varepsilon/2)=
\\
&=v_{mn}+(1-P_2(t_{mn}))\varepsilon/2<v_{mn}+\varepsilon.
  \end{align*}
Estimate the second integrand using the relations~\eqref{lt.2},
\eqref{lt.1}, \eqref{d.1}, \eqref{d.2}, and~\eqref{6a}:
\begin{align*}
&A_1P_1(\xi)+(1-P_1(\xi))\overline{K}_1(\tau^2,y^{2,\varepsilon})=
\\
&=(A_1-\overline{K}_1(\tau^2,y^{2,\varepsilon}))P_1(\xi)+
\overline{K}_1(\tau^2,y^{2,\varepsilon})\le
\\
&\le (A_1-\overline{K}_1(\tau^2,y^{2,\varepsilon}))(P_1(t_{mn})+
\lambda\varepsilon)+\overline{K}_1(\tau^2,y^{2,\varepsilon})\le
\\
&\le A_1P_1(t_{mn})+(1-P_1(t_{mn}))\overline{K}_1
(\tau^2,y^{2,\varepsilon})+(A_1+B_1)\lambda\varepsilon\le
\\
&\le A_1P_1(t_{mn})+(1-P_1(t_{mn}))(v_{m-1,n}+\varepsilon/2)+
\varepsilon/2=
\\
&=v_{mn}+(1-P_2(t_{mn}))\varepsilon/2+\varepsilon/2<v_{mn}+\varepsilon.
  \end{align*}
After integrating we get:
\begin{align*}
&\overline{K}_1(\tau,y^{\varepsilon})< v_{mn}+\varepsilon.
  \end{align*}
\item
Suppose $t\in [t_{mn}+\delta,1]$.
In this case Player II always acts first; so
\begin{align*}
&\overline{K}_1(\tau,y^{\varepsilon})=
\int\limits_{t_{mn}}^{t_{mn}+\delta} (-B_1P_2(\xi)+ (1-P_2(\xi))
\overline{K}_1(\tau^1,y^{1,\varepsilon}) )
d\varphi_{mn}^{\varepsilon}(\xi).
\end{align*}
Estimate the integrand using the relations~\eqref{lt.1}, \eqref{lt.3},
\eqref{6b}:
\begin{align*}
&-B_1P_2(\xi)+ (1-P_2(\xi)\overline{K}_1(\tau^1,y^{1,\varepsilon})=
\\
&=\overline{K}_1(\tau^1,y^{1,\varepsilon})-
(\overline{K}_1(\tau^1,y^{1,\varepsilon})+B_1)P_2(\xi)\le
\\
&\le\overline{K}_1(\tau^1,y^{1,\varepsilon})-
(\overline{K}_1(\tau^1,y^{1,\varepsilon})+B_1)P_2(t_{mn})=
\\
&=-B_1P_2(t_{mn})+
(1-P_2(t_{mn}))\overline{K}_1(\tau^1,y^{1,\varepsilon})\le
\\
&\le-B_1P_2(t_{mn})+
(1-P_2(t_{mn}))(v_{m,n-1}+\varepsilon/2) < v_{mn}+\varepsilon.
  \end{align*}
After integrating we have:
\begin{align*}
&\overline{K}_1(\tau,y^{\varepsilon})< v_{mn}+\varepsilon.
\end{align*}
\end{enumerate}
Hence the inequality~(\ref {t1.5}) is proved.

Now let us turn to the inequality~\eqref{t1.6}.
Let $\eta$ be an arbitrary pure strategy of Player II and
$u=\eta_n$ be the first action time in $\eta$.

There are three cases to be considered.

\begin{enumerate}
\item
Suppose $u\in [0,t_{mn}]$.
In this case Player II acts first; so
\begin{align*}
&\overline{K}_1(x^{\varepsilon},\eta)=
-B_1P_2(u)+(1-P_2(u))\overline{K}_1(x^{2,\varepsilon},\eta^2).
  \end{align*}
By the relations~\eqref{lt.1}, \eqref{lt.3}, and \eqref{6b},
we have:
\begin{align*}
\overline{K}_1(x^{\varepsilon},\eta)=
&
\overline{K}_1(x^{2,\varepsilon},\eta^2)-
(B_1+\overline{K}_1(x^{2,\varepsilon},\eta^2))P_2(u)\ge
\\
&\ge\overline{K}_1(x^{2,\varepsilon},\eta^2)+
(B_1+\overline{K}_1(x^{2,\varepsilon},\eta^2))P_2(t_{mn})=
\\
&=-B_1P_2(t_{mn})+(1-P_2(t_{mn}))\overline{K}_1(x^{2,\varepsilon},\eta^2)
\ge
\\
&\ge -B_1P_2(t_{mn})+(1-P_2(t_{mn}))(v_{m,n-1}-\varepsilon/2)=
\\
&=v_{mn}-(1-P_1(t_{mn}))\varepsilon/2\ge v_{mn}-\varepsilon.
  \end{align*}
\item
Suppose $u\in (t_{mn},t_{mn}+\delta)$.
In this case either player  may act first; so
\begin{align*}
\overline{K}_1(x^{\varepsilon},\eta)&=
\int_{t_{mn}}^u
(A_1P_1(\xi)+(1-P_1(\xi))\overline{K}_1(x^{2,\varepsilon},\eta^2) )
d\varphi_{mn}^{\varepsilon}(\xi)+
\\
&+\int_u^{t_{mn}+\delta}
(-B_1P_2(\xi)+(1-P_2(\xi))\overline{K}_1(x^{1,\varepsilon},\eta^1) )
d\varphi_{mn}^{\varepsilon}(\xi).
  \end{align*}
Estimate the first integrand using the relations~\eqref{lt.2},
\eqref{lt.6}, and~\eqref{6a}:
\begin{align*}
&A_1P_1(\xi)+(1-P_1(\xi))\overline{K}_1(x^{2,\varepsilon},\eta^2)=
\\
&=(A_1-\overline{K}_1(x^{2,\varepsilon},\eta^2))P_1(\xi)+
\overline{K}_1(x^{2,\varepsilon},\eta^2)\ge
\\
&\ge (A_1-\overline{K}_1(x^{2,\varepsilon},\eta^2))P_1(t_{mn})+
\overline{K}_1(x^{2,\varepsilon},\eta^2)=
\\
&= A_1P_1(t_{mn})+(1-P_1(t_{mn}))\overline{K}_1
(x^{2,\varepsilon},\eta^2)\ge
\\
&\ge A_1P_1(t_{mn})+(1-P_1(t_{mn}))(v_{m-1,n}-\varepsilon/2)=
\\
&=v_{mn}-(1-P_2(t_{mn}))\varepsilon> v_{mn}-\varepsilon.
  \end{align*}
Estimate the second integrand using the relations~\eqref{lt.1},
\eqref{lt.2}, \eqref{d.1}, \eqref{d.2}, and~\eqref{6b}:
\begin{align*}
&-B_1P_2(\xi)+(1-P_2(\xi))\overline{K}_1(x^{1,\varepsilon},\eta^1)=
\\
&=\overline{K}_1(x^{1,\varepsilon},\eta^1)-
P_2(\xi)(\overline{K}_1(x^{1,\varepsilon}),\eta^1)+B_1)\ge
\\
&\ge\overline{K}_1(\tau^1,y^{1,\varepsilon})-
(P_2(t_{mn})+\lambda\varepsilon)
(\overline{K}_1(x^{1,\varepsilon},\eta^1)+B_1)\ge
\\
&\ge -B_1P_2(t_{mn})+(1-P_2(t_{mn})) \overline{K}_1
(x^{1,\varepsilon},\eta^1) - (A_1+B_1)\lambda\varepsilon\ge
\\
&\ge -B_1P_2(t_{mn})+(1-P_2(t_{mn}))
(v_{m,n-1}-\varepsilon/2)-\varepsilon/2)=
\\
&=v_{mn}-(1-P_2(t_{mn}))\varepsilon/2-\varepsilon/2> v_{mn}-\varepsilon.
\end{align*}
After integrating we get:
\begin{align*}
&\overline{K}_1(\tau,y^{\varepsilon})> v_{mn}-\varepsilon.
\end{align*}
\item
Suppose $u\in [t_{mn}+\delta,1]$.
In this case Player I acts first; so
\begin{align*}
&\overline{K}_1(x^{\varepsilon},\eta)=
\int\limits_{t_{mn}}^{t_{mn}+\delta}
(A_1P_1(\xi)+(1-P_1(\xi))\overline{K}_1(x^{1,\varepsilon},\eta^1) )
d\varphi_{mn}^{\varepsilon}(\xi).
  \end{align*}
By the relations~\eqref{lt.2}, \eqref{lt.3}, and \eqref{6a},
we have:
\begin{align*}
\overline{K}_1(x^{\varepsilon},\eta)&
=\int\limits_{t_{mn}}^{t_{mn}+\delta}
((A_1-\overline{K}_1(x^{1,\varepsilon},\eta^1))P_1(\xi)+
\overline{K}_1(x^{1,\varepsilon},\eta^1))
d\varphi_{mn}^{\varepsilon}(\xi)\ge
\\
&\ge\int\limits_{t_{mn}}^{t_{mn}+\delta}
((A_1-\overline{K}_1(x^{1,\varepsilon},\eta^1))P_1(t_{mn})+
\overline{K}_1(x^{1,\varepsilon},\eta^1) )
d\varphi_{mn}^{\varepsilon}(\xi)=
\\
&=A_1P_1(t_{mn})+(1-P_1(t_{mn}))\overline{K}_1
(x^{1,\varepsilon},\eta^1)\ge
\\
&\ge A_1P_1(t_{mn})+(1-P_1(t_{mn}))(v_{m-1,n}-\varepsilon/2)=
\\
&= v_{mn}-(1-P_1(t_{mn})\varepsilon/2\ge v_{mn}-\varepsilon/2.
  \end{align*}
\end{enumerate}
Hence the inequality (\ref {t1.6}) is proved.
\end{proof}

Let us continue to prove the theorem. Choose $\varepsilon>0$.
By Lemma~\ref{l.t1}, we can find a strategy $y^{\varepsilon}$,
satisfying (\ref {t1.5}) for any pure strategy $\tau$ of Player I.

From the inequality (\ref {t1.5}) it follows that
\begin{align}
\overline{K}_1(x, y^{\varepsilon})< v_{mn}^1+\varepsilon
\label {t1.7}
  \end{align}
for any mixed strategy $x$ of Player I.
In particular for $x=x^{\varepsilon}$ we have:
\begin{align}
\overline{K}_1(x^{\varepsilon}, y^{\varepsilon})<
v_{mn}^1+\varepsilon\implies v_{mn}^1>\overline{K}_1
(x^{\varepsilon}, y^{\varepsilon})-\varepsilon.
\label {t1.8}
  \end{align}
From the inequality (\ref {t1.6}) it follows that
\begin{align}
\overline{K}_1(x^{\varepsilon},y)> v_{mn}^1-\varepsilon
\label {t1.9}
  \end{align}
for any mixed strategy $y$ of Player II.
In particular for $y=y^{\varepsilon}$ we have:
\begin{align}
\overline{K}_1(x^{\varepsilon}, y^{\varepsilon})>
v_{mn}^1-\varepsilon\implies v_{mn}^1<\overline{K}_1
(x^{\varepsilon}, y^{\varepsilon})+\varepsilon.
\label {t1.10}
  \end{align}
Taking in account the relations (\ref {t1.7}), (\ref {t1.10})
we obtain
\begin{align*}
\overline{K}_1(x, y^{\varepsilon})< v_{mn}^1+\varepsilon<
\overline{K}_1(x^{\varepsilon}, y^{\varepsilon})+\varepsilon
  \end{align*}
for any mixed strategy $x$ of Player I.
Hence inequality (\ref {t1.1}) is proved.
 The relation (\ref {t1.3}) follows from (\ref {t1.8}) and
(\ref {t1.10}).

Due to  the symmetry of the setting the relations (\ref {t1.2}),
(\ref {t1.4}) are proved by the same arguments that
(\ref {t1.1}), (\ref {t1.3}).
\end{proof}

\begin{theorem}   \label{t.1a}
The maxmin value of the game $\overline{\Gamma}_{mn}(P,A,B)$
coincides with the equilibrium value $v_{mn}$.
The $\varepsilon$-equilibrium strategies $x^{\varepsilon}$,
$y^{\varepsilon}$ are $\varepsilon$-maxmin strategies.
\end{theorem}

\begin{proof}
Let us use the properties of the zero-sum duel
$\overline{\Gamma}_{mn}(P,C,\overline C)$, where $C_1=A_1$, $C_2=B_1$,
$\overline C=(C_2,\;C_1)$.
In an antagonistic game the equilibrium value coincides with
the maxmin value and $\varepsilon$-equilibrium strategies are
$\varepsilon$-maxmin strategies.
As the equilibrium value of $\overline{\Gamma}_{mn}(P,C,\overline C)$
is equal to $v_{mn}^1$~\cite{fox}, so  the maxmin value of this game
is also equal to $v_{mn}^1$
and the $\varepsilon$-equilibrium strategy
$x^{\varepsilon}$ of Player I is his $\varepsilon$-maxmin strategy.
Thus the maxmin value of the game $\overline{\Gamma}_{mn}(P,A,B)$
for Player I is equal to $v_{mn}^1$, and $x^{\varepsilon}$
is his $\varepsilon$-maxmin strategy.

The statement of the theorem for Player II follows from the symmetry
of the setting.
\end{proof}

\section{Pareto-optimal plays}

The pair $p=(\tau, \eta)$ of vectors of action moments realized
during the game is called a {\it play}. Let us denote the set of all
plays of a noisy duel $\overline{\Gamma}_{mn}(P,A,B)$ by ${\cal P}$.
Note that ${\cal P}$ is a subset of the set of all situations
of the corresponding silent duel (with the same effectivess functions,
resources, and vectors of profit and loss).
Namely, ${\cal P}$ includes exactly those situations of the silent duel
in which after one of the players has used all of his resource,
the other one postpones his action until the moment $t=1$.
A play $p^1\in {\cal P}$ is called {\it Pareto-optimal\/}
if there exist no $p^2 \in {\cal P}$ such that $p^2\succ p^1$.
Plays $p^1$ and $p^2$ are called {\it incomparable\/} if
$p^1\nsucc p^2$ and $p^2\nsucc p^1$.
Plays $p^1=(\tau^1, \eta^1)$ and $p^2=(\tau^2, \eta^2)$ are called
{\it equivalent\/} if $K_j(\tau^1, \eta^1)=K_j(\tau^2, \eta^2)$ ($j=1,2$).

A play $p=(\tau,\eta)\in {\cal P}$ is called a {\it $T$-play} if
for any $k$, $l$ $(1\leqslant k\leqslant m$; $1\leqslant l\leqslant n)$
at least one of the following two equations holds
$$\tau_k=t_{kl};\quad \eta_l=t_{kl},$$
where $k$, $l$ are the current resources of the players.

\begin{lemma} \label{l.1}
Let $p^T=(\tau^T,\eta^T)$ be an arbitrary $T$-play with noncoinciding
action moments of the players, i.~e., $\tau_k\neq \eta_l$ for all
$k\; (1\leqslant k\leqslant m)$, $l\;(1\leqslant l\leqslant n)$.
Then
\begin{align}
&K_1(\tau^T,\eta^T)=v_{mn}^1;\quad
K_2(\tau^T,\eta^T)=v_{mn}^2.
\label {l3.4}
  \end{align}
\end{lemma}

\begin{proof}
Let us prove Lemma by induction on the number of action moments
of the players.
For $n=0$, $m>0$ or $m=0$, $n>0$ the assertion is true as
$v_{m0}^1=A_1$; $v_{m0}^2=-B_2$, $v_{0n}^1=-B_1$; $v_{0n}^2=A_2$.
 Assume that
the statement is true for all pairs $(k,l)$ such that $k\leqslant m$;
 $l\leqslant n$; $k+l<m+n$ and prove it for $(k,l)=(m,n)$.
Suppose that Player I acts at the moment $t=t_{mn}$.
Then using the recursive formula~\eqref{4a} for $\tau_m<\eta_n$ and
the inductive assumptions we obtain
\begin{align*}
K_1(\tau,\eta)&=
A_1P_1(t_{mn})-\left(1-P_1(t_{mn})\right)
\left(A_1-(A_1+B_1)\prod_{i=1}^{m-1}\left(1-P_1(t_{in})\right)\right)=
\\
&=A_1-(A_1+B_1)\prod_{i=1}^m\left(1-P_1(t_{in})\right)=v_{mn}^1.
  \end{align*}
Suppose that Player II acts at the moment of time~$t_{mn}$.
Then using the formula~\eqref{4a} for $\tau_m>\eta_n$
and the inductive assumptions we obtain
\begin{align*}
K_1(\tau,\eta)&=
-B_1P_2(t_{mn})-\left(1-P_2(t_{mn})\right)
\left((A_1+B_1)\prod_{j=1}^{n-1}\left(1-P_2(t_{mj})\right)-B_1\right)=
\\
&=-B_1+(A_1+B_1)\prod_{j=1}^n\left(1-P_2(t_{mj})\right)=v_{mn}^1.
  \end{align*}
The second player's payoff fuction is considered in the analogous way.
\end{proof}

\begin{lemma}          \label {l.4}
Let the plays $p^1=(\tau^1,\eta^1)$, $p^2=(\tau^2,\eta^2)\in{\cal P}$
satisfy the conditions
 \begin{align}
 & \eta^1_i=\eta^2_i, \ \tau^1_i=\tau^2_i \quad\text{ for } i\geqslant 2;
\label {l4.1}\\
& \tau^1_1=t_{11};\quad \eta^1_2<t_{11};\quad \eta^1_1=1;\quad
 \tau^2_1=t_{11};\quad  \eta^2_1=t_{11}.
\label {l4.3}
 \end{align}
Then we have
(1) if $A\succ B$, then $p^1\succ p^2$;
(2) if $B\succ A$, then $p^2\succ p^1$.
\end{lemma}

\begin{proof}
Express the payoff function of the plays $p^1$ and $p^2$
in the following way:
\begin{align*}
K_1(\tau^1,\eta^1)&=K_1(\tau_m,\dots\tau_2;\eta_n,\dots\eta_2)+
\prod_{i=2}^m\left(1-P_1(\tau_{i})\right)\times
\\
&\times\prod_{j=2}^n\left(1-P_2(\eta_{j})\right)
\left(A_1P_1(t_{11})-B_1\left(1-P_1(t_{11})\right)\right);
\\
K_1(\tau^2,\eta^2)&=K_1(\tau_m,\dots\tau_2;\eta_n,\dots\eta_2)+
\prod_{i=2}^m\left(1-P_1(\tau_{i})\right)\times
\\
&\times\prod_{j=2}^n\left(1-P_2(\eta_{j})\right)
\left(A_1P_1(t_{11})\left(1-P_2(t_{11})\right)-
B_1\left(1-P_1(t_{11})\right)P_2(t_{11})\right).
\end{align*}
Consider the difference
$K_1(\tau^1,\eta^1)-K_1(\tau^2,\eta^2)$. Taking in account that by
 Lemma~\ref{l.0}
$$P_1(t_{11})+P_2(t_{11})=1,$$
we have:
\begin{align}
K_1(\tau^1,\eta^1)&-K_1(\tau^2,\eta^2)=
\nonumber
\\
&= (A_1-B_1)P_1(t_{11})P_2(t_{11})
\prod_{i=2}^m\left(1-P_1(\tau_{i})\right)
\prod_{j=2}^n\left(1-P_2(\eta_{j})\right).
\label {l3.1}
  \end{align}
Analogously,
\begin{align}
K_2(\tau^1,\eta^1)&-K_2(\tau^2,\eta^2)=
\nonumber
\\
&=(A_2-B_2)P_1(t_{11})P_2(t_{11})
\prod_{i=2}^m\left(1-P_1(\tau_{i})\right)
\prod_{j=2}^n\left(1-P_2(\eta_{j})\right).
\label {l3.2}
  \end{align}
The statement of Lemma follows from~\eqref{l3.1}, \eqref{l3.2}.
 \end{proof}

 A duel may have one of four alternative results:
$H_0$ --- no one of players achieves success;
$H_1$ --- Player I achieves success;
$H_2$ --- Player II achieves success;
$H_{3}$ --- both players achieve success simultaneously.
We denote the probability of the result $H_i$ in play $(\tau,\eta)$ by
 $Q_i(\tau,\eta)$ ($i=0,1,2,3$). Then:
\begin{align}
\sum_{i=0}^{3}Q_i(\tau,\eta)=1.
\label{eqno3}
  \end{align}
Success of both players is possible only if they act simultaneously,
because if one of the players achieves success, the game stops.
As the payoff function $K_j(\tau; \eta)$ is the mathematical
expectation of profit, so
\begin{align}
&K_1(\tau; \eta)= A_1Q_1(\tau,\eta)-B_1Q_2(\tau,\eta);
\label{eqno4}
\\
&K_2(\tau; \eta)= A_2Q_2(\tau,\eta)-B_2Q_1(\tau,\eta).
\label{eqno4a}
  \end{align}
By~\eqref{eqno3},  we have:
\begin{align}
&K_1(\alpha,\beta)=A_1-(A_1+B_1)Q_2(\alpha,\beta)-
A_1(Q_0(\alpha,\beta)+Q_3(\alpha,\beta));
\label{eqno5}
\\
&K_2(\alpha; \beta)=-B_2+(A_2+B_2)Q_2(\alpha,\beta)+
B_2(Q_0(\alpha,\beta)+Q_3(\alpha,\beta)).
\label{eqno6}
  \end{align}

 Denote by ${\cal P}'\subset {\cal P}$ the set of plays with
noncoinsiding action times in which $\tau_1=1$ or $\eta_1=1$.
 The complement ${\cal P}\setminus {\cal P}'$ includes plays in
which the players use their last units of resource simultaneously.
\begin{lemma}
          \label {l.5}
Let  $p'=(\tau',\eta') \in {\cal P}'$ and
 $p=(\tau,\eta)\in {\cal P}\setminus {\cal P}'$.
In this case

(1) if  $p'\succ p$, then $A_1A_2>B_1B_2$;\quad
(2) if $p\succ p'$, then $A_1A_2<B_1B_2$.
\vskip 1mm
\end{lemma}

\begin{proof}
Let $(\tau',\eta')\in {\cal P}'$, so
the relations~\eqref{eqno5}, \eqref{eqno6} reduce to the form:
  \begin{align}
&K_1(\tau', \eta')=
A_1-(A_1+B_1)Q_2(\tau',\eta');
\label{eqno7}
\\
&K_2(\tau', \eta')=
-B_2+(A_2+B_2)Q_2(\tau',\eta').
\label{eqno7a}
  \end{align}
Consider the differences $\Delta_j=K_j(\tau,\eta)-K_j(\tau',\eta')$.
Applying the formulas~\eqref{eqno5}, \eqref{eqno6}, \eqref{eqno7},
\eqref{eqno7a}
 we get:
  \begin{align}
&\Delta_1=(A_1+B_1)I(\tau',\eta';\tau,\eta)-A_1Q(\tau,\eta);
\label{eqno8}
\\
&\Delta_2=-(A_2+B_2)I(\tau',\eta';\tau,\eta)+B_2Q(\tau,\eta).
\label{eqno9}
  \end{align}
where
  \begin{align*}
&I(\tau',\eta';\tau,\eta)=
Q_2(\tau',\eta')-Q_2(\tau,\eta);
\\
&Q(\tau',\eta')=Q_0(\tau',\eta')+Q_3(\tau',\eta').
  \end{align*}
  Let us prove the first assertion of Lemma.
 If $p'\succ p$, then  $\Delta_1\le 0$, $\Delta_2\le 0$ and
at least one of this inequalities is strict.
 Since $A_j+B_j>0$ ($j=1,2$), it follows from~\eqref{eqno8},
\eqref{eqno9} that
  \begin{align}
&I(\tau',\eta';\tau,\eta)\le\frac{A_1}{A_1+B_1}Q(\tau,\eta);
\label{eqno11}
\\
&I(\tau',\eta';\tau,\eta)\ge\frac{B_2}{A_2+B_2}Q(\tau,\eta).
\label{eqno11a}
  \end{align}
 Since $Q(\tau,\eta)\neq 0$ for
$p=(\tau,\eta)\in {\cal P}\setminus {\cal P}'$,
combining the inequalities~\eqref{eqno11}, \eqref{eqno11a}
and taking in account that
at least one of them is strict, we get:
    $$\frac{B_2}{A_2+B_2}<\frac{A_1}{A_1+B_1}\implies A_1A_2>B_1B_2.$$

  The second statement of the lemma is proved in the analogous way.
 \end{proof}

\begin{lemma}
      \label {l.6}
Let  $p_1, p_2 \in {\cal P}'$.
Then $p_1$ and $p_2$ are Pareto-incomparable.
\end{lemma}

\begin{proof}
If $p=(\tau,\eta)\in {\cal P}'$, then from
the relations~\eqref{eqno7}, \eqref{eqno7a} it follows that:
  \begin{align}
K_2(\alpha,\beta)=\frac{A_1A_2-B_1B_2}{A_1+B_1}-
\frac{A_2+B_2}{A_1+B_1}K_1(\alpha,\beta).
\label{eqno12}
  \end{align}
The statement of Lemma immediately follows from \eqref{eqno12}.
 \end{proof}

\begin{theorem}
\label{t.6}
If in the game $\Gamma_{mn}(P,A,B)$ the coefficients of profit
and loss of the players satisfy the condition
 \begin{align}
A_1A_2=B_1B_2,
\label {P2}
 \end{align}
then the game is quasi-antagonistic.
\end{theorem}

In the paper~\cite{pareto} a similar theorem about a sufficient
condition of quasi-antagonisticity was also proven for continuous
and discrete duels.

\begin{proof}
Note that if the relation~\eqref{P2} holds,
then there exists a number $\lambda>0$ such that
 \begin{align}
K_1(\tau,\eta)=-\lambda K_2(\tau,\eta).
\label{eqno14}
 \end{align}
Indeed, suppose
$A_1A_2>0$. Then accoding to~\eqref{P2}
$B_1B_2>0$ and $A_1/B_2=B_1/A_2$.
In this case from the formulas~\eqref{eqno4}, \eqref{eqno4a}
one can see that~\eqref{eqno14} holds with $\lambda=A_1/B_2$.
Suppose that one of the numbers $A_1$ or $A_2$ is equal to zero.
If $A_1=0$, then by~\eqref{I1} we have $B_1\neq 0$, and
from~\eqref{P2} it follows $B_2=0$, $A_2\neq 0$.
Then in view of the relations~\eqref{eqno4}, \eqref{eqno4a} we conclude
that the equation~\eqref{eqno14} holds with $\lambda=B_1/A_2$.
The case $A_2=0$ is considered in the similar way.
So it is proved that there exists $\lambda>0$ for which
the equation~\eqref{eqno14} holds.
It follows from~\eqref{eqno14} that the game is quasi-antagonistic.
 \end{proof}

\begin{theorem}   \label{t.4}
If the coefficients of profit and loss of the players
are related by the inequality
$A_1A_2\geqslant B_1B_2$,
then $T$-plays with noncoinciding action moments of the players
are Pareto-optimal.
\end{theorem}

\begin{proof}
If $A_1A_2=B_1B_2$, then by Theorem~\ref{t.6}
all plays are Pareto-optimal.
In any $T$-play with noncoinciding action moments
the last action of a player happens at $t=1$, i.~e.,
 $\max\{\tau_{1},\eta_{1}\}=1$.
Suppose $A_1A_2>B_1B_2$. We will show that $p\nsucc p'$
for any plays $p\in {\cal P}$, $p'\in {\cal P}'$.
There are two cases to be considered.
\begin{enumerate}
\item
$p\in {\cal P}'$. Then by Lemma~\ref{l.6} $p\nsucc p'$.
\item
$p\in {\cal P}\setminus {\cal P}'$.
By Lemma~\ref {l.5}, if $p\succ p'$ then $A_1A_2<B_1B_2$ in
contradiction to the assumption of the theorem.
 Therefore $p\nsucc p'$.
\end{enumerate}
 \end{proof}

\begin{theorem}   \label{t.4a}
If $B\succ A$, then $T$-plays with noncoinciding action moments
of the players are not Pareto-optimal.
\end{theorem}

 \begin{proof}
Since according to Lemma~\ref{l.1} all the $T$-plays with noncoinciding
action moments of the players are equivalent, we will give the proof
for one $T$-play only, namely for the play
 $p^1=(\tau, \eta)$ in which
 \begin{align*}
&\tau_i=t_{i1} \text { for }i=1,\dots, m;\\
& \eta_j=t_{mj} \text { for }j=2,\dots, n;\quad \eta_1=1.
 \end{align*}
Set
$$p^2=(\tau, \eta^2), \text{ where }\ \eta^2=(t_{11},\eta_2,\dots,\eta_n).$$
The plays $p^1, p^2$ satisfy the conditions of Lemma~\ref{l.4},
therefore $$B\succ A\implies p^2\succ p^1.$$
Thus the $T$-play $p^1$ and all the other $T$-plays with noncoinciding
action moments of the players are not Pareto-optimal.
 \end{proof}

\begin{theorem}
 \label{t.5}
If the duel
$\Gamma_{mn}(P,A,B)$ is quasi-antagonistic, then
one of the following two conditions holds:
 \begin{align}
(1)\; (A_1- B_1)(A_2-B_2)<0;\qquad
(2)\; A_j=B_j; \;j=1,2.
\label {P1}
 \end{align}
\end{theorem}

 \begin{proof}
Assume that the plays $p^1, p^2$ satisfy
the conditions~\eqref{l4.1}, \eqref{l4.3}.
Suppose that no one of the conditions~(\ref{P1}) holds.
Then the following two cases are possible.
\begin{enumerate}
\item
 $A\succ B$.
Then according to the statement~(1) of Lemma~\ref{l.4} one has
$p^1\succ p^2$. Therefore the play $p^2$ is not Pareto-optimal,
and so the game is not quasi-antagonistic.
\item
 $B\succ A$.
Then according to the statement~(2) of Lemma~\ref{l.4} one has
$p^2\succ p^1$. Therefore the play $p^1$ is not Pareto-optimal,
and so the game is not quasi-antagonistic.
\end{enumerate}
 \end{proof}

\end{document}